\documentclass[letterpaper]{article} 
\usepackage{aaai25}
\usepackage{times}  
\usepackage{helvet}  
\usepackage{courier}  
\usepackage[hyphens]{url}  
\usepackage{graphicx} 
\usepackage{natbib}  
\usepackage{caption} 
\usepackage{algorithm}

\usepackage{graphicx}
\usepackage{xcolor}
\usepackage{subcaption}
\usepackage{cite}
\usepackage{amsmath}
\usepackage{enumitem}
\usepackage{subcaption}
\usepackage[utf8]{inputenc}
\usepackage{algpseudocode}
\usepackage{booktabs}
\usepackage{multirow}
\usepackage{tabularx}
\usepackage{parskip}
\usepackage{caption}
\usepackage{soul}
\usepackage{nicematrix}
\usepackage{booktabs}
\usepackage{makecell}  
\usepackage{comment}
\usepackage{subcaption}
\usepackage{amssymb}
\usepackage{mathtools}

\usepackage{newfloat}
\usepackage{listings}
\DeclareCaptionStyle{ruled}{labelfont=normalfont,labelsep=colon,strut=off} 
\floatstyle{ruled}
\newfloat{listing}{tb}{lst}{}
\floatname{listing}{Listing}
\title{Efficient Primal Heuristics for Mixed Binary Quadratic Programs \\
Using Suboptimal Rounding Guidance}
\author{
    Weimin Huang\textsuperscript{\rm 1},
    Natalie M. Isenberg\textsuperscript{\rm 2},
    Ján Drgoňa\textsuperscript{\rm 3},
    Draguna L Vrabie\textsuperscript{\rm 2},
    Bistra Dilkina\textsuperscript{\rm 1}
}
\affiliations{
    \textsuperscript{\rm 1}University of Southern California, USA\\
    \textsuperscript{\rm 2}Pacific Northwest National Laboratory, USA\\
    \textsuperscript{\rm 3}Johns Hopkins University, USA\\

    weiminhu@usc.edu, natalie.isenberg@pnnl.gov, jdrgona1@jh.edu, draguna.vrabie@pnnl.gov, dilkina@usc.edu
}

\usepackage{bibentry}

\begin{document}

\maketitle

\begin{abstract}
Mixed Binary Quadratic Programs (MBQPs) are a class of NP-hard problems that arise in a wide range of applications, including finance, machine learning, and chemical and energy systems. Large-scale MBQPs are challenging to solve with exact algorithms due to the combinatorial search space and nonlinearity. 
Primal heuristics have been developed to quickly identify high-quality solutions to challenging combinatorial optimization problems. In this paper, we propose an extension for two well-established rounding-based primal heuristics, RENS and Undercover. Instead of using the optimal solution to a relaxation for variable rounding and search as in RENS, we use a suboptimal relaxation solution of the MBQP as the basis for rounding and guidance for searching over a restricted subproblem where a certain percentage of binary variables are free. We apply a similar idea to the Undercover heuristic that fixes a variable cover to the rounded relaxation values. Instead, we relax a subset of the cover variables based on the suboptimal relaxation and search over a larger restricted subproblem. We evaluate our proposed methods on synthetic MBQP benchmarks and real-world wind farm layout optimization problem instances. The results show that our proposed heuristics identify high-quality solutions within a small time limit and significantly reduce the primal gap and primal integral compared to RENS, Undercover, and solvers with additional primal heuristics integrated inside Branch-and-Bound. 
\end{abstract}

%

\section{Introduction}



Mixed Binary Quadratic Programs (MBQPs) are mathematical optimization problems with binary variables that contain quadratic terms in the objective function. MBQPs are classic problems in combinatorial optimization that can capture real-world applications including finance~\cite{parpas2006global}, machine learning~\cite{bertsimas2009algorithm}, and chemical~\cite{misener2013glomiqo} and energy systems~\cite{turner2014new}. Many important combinatorial optimization problems can be formulated as MBQPs, including the quadratic assignment problem~\cite{loiola2007survey}, the stable set problem~\cite{rebennack2024stable}, and the vertex coloring problem~\cite{kochenberger2005unconstrained}. 

MBQPs are NP-hard in general \cite{pia2017mixed} and are particularly challenging to solve due to the combinatorial search space coupled with nonlinearities. The Branch-and-Bound (BnB) algorithm is an exact tree search algorithm for solving MBQPs and more general Mixed-Integer Nonlinear Programming (MINLP) problems. Although there is a large body of work to improve BnB algorithms, large-scale MBQPs (numbering in $10^2 - 10^4$ or more binary variables) are still challenging to solve with exact methods \cite{silva2021quadratic}. 

Correspondingly, a significant body of research has focused on 
\emph{primal heuristics}. Primal heuristics are algorithms designed to quickly identify high-quality feasible solutions for a given optimization problem without optimality guarantees \cite{berthold2014rens}. As far as the authors are aware, there has been limited work in the existing literature on primal heuristics specifically designed for solving MBQPs. Many primal heuristics are initially developed for Mixed Integer Linear Programming (MILP) problems and have been adapted for general MINLPs \cite{berthold2014rens,bonami2009feasibility}. Given the prevalence of MBQP formulations in industrial applications, the challenges in solving MBQPs via exact and heuristic methods \cite{silva2021quadratic}, and the demand for real-time solutions to MBQPs in scientific domains (e.g., model predictive control) \cite{takapoui2020simple}, there is a need for efficient algorithms that can quickly find high-quality solutions to these problems.

In this work, we develop efficient primal heuristics for large-scale general MBQPs based on relaxation, rounding, and search techniques. We propose \emph{Relax-Search} and \emph{Cover-Relax-Search}, which are extensions of two well-established rounding-based primal heuristics. While RENS \cite{berthold2014rens} solves a relaxation of the original problem to optimality, fixes the subset of integer variables that are integral in the relaxation solution and uses a complete solver to search over the resulting subproblem, Relax-Search finds a suboptimal solution to the relaxation in limited time and uses it as guidance for the partial variable fixing. 
We apply a similar idea to the Undercover heuristic \cite{berthold2014undercover} that solves a relaxation, finds a variable \textit{cover} that removes nonlinearities, and fixes the variables in the \textit{cover} to the rounded relaxation values. We propose the Cover-Relax-Search method which computes a suboptimal variable \textit{cover} and a suboptimal relaxation solution in limited time and uses them to fix a subset of the \textit{cover} variables and search over a larger restricted subproblem.




We evaluate the proposed approaches on multiple standard MBQP benchmarks. Additionally, we demonstrate the performance of the proposed methods on a real-world Wind Farm Layout Optimization Problem (WFLOP). We extend an existing MBQP formulation \cite{turner2014new} for WFLOP to account for uncertain wind conditions and construct real-world instances based on wind distribution data in California. Results across all benchmarks -- synthetic and real-world -- show that our proposed approach identifies high-quality feasible solutions in a short time cutoff and significantly reduces the primal integral compared to the RENS
heuristic, the Undercover heuristic, and solvers with additional primal heuristics integrated inside Branch-and-Bound.

\section{Background}
\subsection{Mixed Binary Quadratic Programs}
Mixed Binary Quadratic Programs (MBQPs) are the class of problems of the form
\[
\min  x^T H x + c^T x 
\]
\[
\text{ s.t. } A x \leq b
\]
\[
x_j \in \{0,1\}, \forall j \in B,
\]
where $c \in \mathbb{R}^n$, $b \in  \mathbb{R}^m$, $H \in  \mathbb{R}^{n \times n}$, and $A \in  \mathbb{R}^{m\times n}$. $B \subseteq \{1,...,n\}$ is the set of variables that are restricted to be binary. $H$ is a real symmetric matrix and is not necessarily positive semidefinite, allowing for nonconvex objective functions. MBQPs are NP-hard in general \cite{pia2017mixed}. The Branch-and-Bound (BnB) algorithm is an exact tree search algorithm to solve MILPs, MBQPs, and general MINLPs. 

\subsection{Reformulation and Relaxation of MBQPs}
\paragraph{Nonlinear Programming Relaxation} The Nonlinear Programming (NLP) relaxation of the MBQP is obtained when the integrality constraints on the set of binary variables $B$ are relaxed such that $x_j \in [0,1], \forall j \in B$. The optimal solution to the NLP relaxation is not necessarily integral, but the optimal value of the relaxed problem serves as a lower bound to the MBQP, as the feasible
region of the NLP relaxation problem is larger than the MBQP.

\paragraph{Linear Reformulation and Relaxation} A common approach to solving MBQPs is to use a linearization step that reformulates the nonlinear terms in the objective function into an equivalent linear form by introducing auxiliary variables and constraints \cite{forrester2020computational}. Terms in MBQPs that are products of binary variables $x_i x_j$, $i,j \in B$ are replaced by a new continuous variable $z_{ij}$. 
Based on the Reformulation-Linearization Technique (RLT) \cite{sherali2013reformulation}, the following McCormick (RLT-0) linear constraints are added:
\[
z_{ij} \geq x_i +x_j -1, \forall i<j \]
\[
z_{ij} \leq x_i, \forall i<j
\]
\[
z_{ij} \leq x_j, \forall i<j
\]
Terms $x_i^2, i \in B$ are replaced by $x_i$ as $x_i = x_i^2$ by exploiting the binary property. With this reformulation, the Linear Programming (LP) relaxation of a MBQP can be obtained by relaxing the binary variables to $x_j \in [0,1], \forall j \in B$ in the reformulation problem. The solutions obtained by solving this relaxation also provide a valid lower bound to the true optimal value, but may not be integral feasible.

\subsection{Primal Heuristics}

Primal heuristics are algorithms designed to quickly find high-quality solutions for a given optimization problem without optimality guarantees \cite{berthold2014rens}. Primal heuristics can be used as either a standalone algorithm or as a supplementary procedure that improves the upper bound in BnB. To our knowledge, there has been limited research on primal heuristics specifically designed to solve MBQPs. Most existing heuristics were initially developed for MILPs and later adapted for general MINLPs \cite{berthold2014rens,bonami2009feasibility}. The LP reformulation described in the previous section enables the application of MILP heuristics to MBQPs \cite{berthold2014primal}. 

Many of these MILP and MINLP heuristics are based on rounding and searching over a restricted subproblem. Rounding-based heuristics use the continuous relaxation of the original problem as a reference solution \cite{berthold2014primal}. These heuristics typically involve solving the relaxation and then creating a subproblem by fixing a subset of integer variables by rounding the relaxation values to the nearest integer values \cite{nannicini2012rounding,berthold2014rens,takapoui2020simple}. The proposed methods in this work build upon two well-established rounding-based heuristics, RENS and Undercover, both of which are start heuristics that do not require an initial feasible solution.



\paragraph{RENS}\label{rens} The Relaxation Enforced Neighborhood Search (RENS) algorithm is a rounding-based heuristic that uses an optimal solution $\bar{x}$ of an LP or NLP relaxation as a reference solution and creates a MINLP subproblem to be solved with a complete solver. The generic RENS algorithm is described in Algorithm \ref{alg:rens}. The key idea of the algorithm is to fix variables that take integral relaxation values (Line 4) and optimize over all the possible rounding of the fractional variables by solving a sub-MINLP (Line 9). The MINLP subproblem is constructed by fixing integer variables for which the relaxation is integral and changing the bounds of the remaining integer variables to the two nearest integers. For the case of MBQPs, RENS fixes all variables for which the relaxation take binary values, and the bounds for the remaining binary variables remain unchanged. Beyond the generic implementation shown in Algorithm \ref{alg:rens}, the full-fledged RENS algorithm includes additional components such as probing and conflict analysis to improve primal solutions in the sub-MINLP \cite{berthold2014rens}.

\begin{algorithm}
\caption{Relaxation Enforced Neighborhood Search (RENS)}
\label{alg:rens}
\begin{algorithmic}[1]
\Require MINLP $\mathcal{P}$ with set of integer variables $\mathcal{I}$
\State $\bar{x} \gets$  Compute optimal solution of the relaxation of $\mathcal{P}$ 
\For {$i \in \mathcal{I}$}
    \If {$\bar{x_i}$ is integral}
        \State Fix $x_i = \bar{x_i}$
    \Else
        \State Restrict the bounds of $x_i$ to $\lfloor \bar{x_i} \rfloor \leq x_i \leq \lceil \bar{x_i} \rceil$
    \EndIf
\EndFor
\State Solve the restricted sub-MINLP using a MINLP solver

\end{algorithmic}
\end{algorithm}

\paragraph{Undercover}\label{undercov}
Most MINLP heuristics require solving a smaller subproblem that is of the same class as the original optimization problem \cite{berthold2014rens}. However, MINLP subproblems derived from a larger MINLP are not necessarily easier to solve than the original problem, as both the discrete nature and nonlinear terms contribute to the complexity in the subproblem. The Undercover heuristic creates a subproblem of an easier class compared to the original problem \cite{berthold2014undercover, belotti2013mixed}. Specifically, it identifies the minimal set of variables needed to be fixed to remove the nonlinear terms. This set of variables is called the \textit{cover} set and can be computed as the minimum vertex cover of the Hessian graph. Formally, a Hessian graph for a given MINLP instance is a graph $G=(V,E)$, where $V$ is the set of variables. There is an edge $e=(i,j)$ between variable $i$ and variable $j$ if and only if $i$ and $j$ appear together in a nonlinear term. The minimum vertex cover of $G$ can be obtained by solving the following optimization problem. Let the binary decision variable $\alpha_i$ be 1 if and only if variable $i$ is in the cover, and solve for
\begin{equation}
\min \sum_{i \in V} \alpha_i \; \text{s.t. } \alpha_i +\alpha_j \geq 1 \; \text{for all}\; (i,j) \in E \label{eq:1}
\end{equation}
By definition of the Hessian graph, the minimum vertex cover of $G$ is the minimal set of variables needed to be fixed to obtain an MILP subproblem that does not contain nonlinear terms from the original MINLP. In Undercover, all integer \textit{cover} variables are fixed to the rounded LP or NLP relaxation values. Then, it solves the MILP subproblem in (\ref{eq:1}). Although this MILP subproblem in (1) is also NP-hard, \cite{berthold2014undercover} show that the runtime to compute the cover is often short compared to the runtime for other components of the algorithm, such as solving the relaxation and the subproblem. 
\section{Efficient Primal Heuristics using Suboptimal Rounding Guidance} \label{sec:CRS}
We present two simple yet effective primal heuristics for MBQPs, \emph{Relax-Search} and \emph{Cover-Relax-Search}, inspired by RENS and Undercover. Both proposed heuristics are start heuristics that do not require a feasible solution to the original problem to start with.

\subsection{Relax-Search}
We propose \emph{Relax-Search}, which is an extension of RENS that uses the suboptimal relaxation as guidance for rounding and variable fixing and searches over a modified restricted subproblem to improve over the suboptimal rounding. While the percentage of variables to be fixed in RENS depends on both the original problem and the relaxation solution, we control the percentage of variables fixed using a ratio $p$. 

The Relax-Search algorithm is shown in Algorithm \ref{myalgo}. First, in the \textit{relax} phase, given a MBQP $\mathcal{P}$, Relax-Search computes the LP or NLP relaxation of $\mathcal{P}$ given a time limit $T_r$. The relaxation solution (denoted as $\bar{x}$) can be suboptimal given the time cutoff, compared to RENS which uses the optimal relaxation. Second, in the \textit{search} phase, the algorithm treats the set of binary variables $\mathcal{B}$ as the candidate fixing set $\mathcal{U}$ and selects a subset $\mathcal{U}' \subseteq \mathcal{U}$ to fix in creating the restricted subproblem. The size $k$ of this subset is the number of variables in the candidate set multiplied by a ratio $p$. To decide on the set of variables to be fixed $\mathcal{U}'$, we compute the degree to which the relaxed solution takes integer values $\Delta_i = |\bar{x_i} - 0.5|$ and select the top $p \times |\mathcal{U}|$ variables with the largest $\Delta_i$. For variables $i \in \mathcal{U}'$, we fix the values to be the rounded relaxation $\lfloor\bar{x_i}\rceil$. The variables in the set $\mathcal{B} \setminus \mathcal{U}' $ are free in the sub-MBQP.


\begin{algorithm}[t]
\caption{Relax-Search {\color{blue} (Cover-Relax-Search)}}
\label{myalgo}
\begin{algorithmic}[1]
\Require  A MBQP $\mathcal{P}$ with set of binary variables $\mathcal{B}$, fixing ratio $p$, relaxation time limit $T_r$ {,\color{blue}(cover time limit $T_c$)} 
\State $\bar{x} \gets$  Compute the relaxation of $\mathcal{P}$ given time limit $T_r$

\State Candidate fixing set $\mathcal{U} \gets \mathcal{B}$
{\color{blue}\State Candidate fixing set $\mathcal{U} \gets$ Compute the \textit{cover} variables of the MBQP by solving the cover problem in (\ref{eq:1}) (Time limit: $T_c$)}
\State $\Delta_i \gets |\bar{x_i}-0.5|$ for all $i\in \mathcal{U}$
\State $k \gets |\mathcal{U}| \times p$ 
\State $\mathcal{U}'$  $\gets$ Select $k$ variables from $\mathcal{U}$ greedily with the largest $\Delta_i$ 
\For {$i \in \mathcal{U}'$}
    \State Fix $x_i = \lfloor\bar{x_i}\rceil $ by rounding to the nearest integer
\EndFor
\State Solve the sub-MBQP using a solver

\end{algorithmic}
\end{algorithm}

\subsection{Cover-Relax-Search}
We apply a similar idea to extend the Undercover heuristic and propose \emph{Cover-Relax-Search}. The relevance to \emph{Relax-Search} is shown in Algorithm ~\ref{myalgo}. In the \textit{cover} phase, we compute the \textit{cover} set variables $C$ for the given MBQP $\mathcal{P}$ by solving optimization problem (1), as in the Undercover heuristic. In the \textit{relax} phase, similar to Relax-Search, we compute the LP or NLP relaxation of $\mathcal{P}$ under a time limit $T_r$, which can be suboptimal. In the \textit{search} phase, Cover-Relax-Search treats the \textit{cover} set $\mathcal{C}$ as the candidate fixing set $\mathcal{U}$ and fixes the top $p \times |\mathcal{U}|$ variables for which the suboptimal relaxation values are closer to integrality. Compared to the original Undercover that fixes all variables in $\mathcal{C}$, we allow a fraction of nonlinear terms in the subproblem and create a larger search space. The original Undercover heuristic follows the "fast fail" strategy, which means that it can compute a solution for the subproblem efficiently in the case when the rounding is feasible and does not consume much running time in the case when the resulting subproblem is infeasible \cite{berthold2014undercover}. Our extension improves the success rate of Undercover as a standalone heuristic by controlling the portion of variables in $C$ that are fixed.

\section{Experiments and Analysis} \label{sec:evaluation}
We divide the empirical evaluation into two main parts. In this section, we evaluate our methods on standard synthetic MBQP benchmarks. In the next section, we examine the efficacy of our method by applying it to a large-scale, real-world wind farm layout optimization problem.

\subsection{Setup}

\subsubsection{Instances Generation} We evaluate the proposed methods on three NP-hard MBQP benchmarks, which are the Cardinality-constrained Binary Quadratic Programs (CBQP)  \cite{zheng2012successive}, Cardinality-constrained Quadratic Knapsack Problem (CQKP) \cite{letocart2014efficient}, and the Quadratic Multidimensional Knapsack Problem (QMKP) \cite{forrester2020computational}. We generate 100 small and 100 large test instances for each problem class. Small and large instances contain 500 and 1000 binary variables, respectively (see Table \ref{tab-stas}). We follow the approach described in \cite{forrester2020computational} to generate random $H$ matrices for general nonconvex MBQPs. 

\begin{table} [!htbp]
\centering

\scalebox{0.9}{
\begin{tabular}{|r|r|r|r|r|}
 \hline
Benchmark & \# Var & \# Cons & Quad. den. & PD Gap\\
 \hline
 CBQP-small & 500 &  1& 0.1& 605.42\% \\
 \hline
 CBQP-large  & 1000 &  1 &  0.1 & 713.51 \%\\
  \hline
CQKP-small &500 &  2 &  0.1& 476.45\%\\
 \hline
 CQKP-large & 1000 &  2 & 0.1&846.66\%\\
  \hline
QMKP-small  & 500 &  50 &   0.1&370.59\% \\
 \hline
QMKP-large &1000 &  50&  0.1&  1022.49\%\\
 \hline
\end{tabular}}
\caption{Instance statistics (100 instances for each benchmark). Number of binary variables (\# Var), number of constraints (\# Cons), quadratic term density (Quad. den,), and Primal-Dual Gap at 30 minutes (PD Gap) with SCIP.}
\label{tab-stas}
\end{table}
\subsubsection{Evaluation Metrics} We use the following metrics to evaluate the effectiveness
of different methods: (1) The \textit{primal gap} \cite{berthold2013measuring} is the normalized difference between the objective value $v$ of the MBQP and a best known objective value $v^*$. In the case where the optimal solution to the MBQP is unknown, the best known solution $v^*$ is the best objective value found in all methods tested given the time limit. When a method has found a solution with objective value $v$, the primal gap is defined as
\[
\text{Primal Gap} =
\begin{cases} 
    0, & \text{if } |v|=|v^*|=0, \\
    1, & \text{if } vv^* <0, \\
    \frac{|v-v^*|}{max(|v|,|v^*|)}, & \text{else}.
\end{cases}
\]
When no feasible solution is found, the primal gap is defined to be 1. (2) The \textit{primal integral} \cite{berthold2013measuring} is the integral of the primal gap over time, which has been widely used in benchmarking primal heuristics in both MILPs and MINLPs. It captures the solution quality and the speed at which better solutions are found. It also serves as a measure of convergence towards the optimal solution (or the best-known solution when the optimal solution is unknown). 

\subsubsection{Baselines}
We compare our methods with the following baselines:
\begin{itemize}
    \item \textbf{SCIP with a portfolio of primal heuristics integrated} \cite{BestuzhevaEtal2021OO}. SCIP uses BnB  as its main algorithm and contains a portfolio of primal heuristics for MINLPs to find good feasible solutions, including Feasibility Pumps, Large Neighbourhood Search, Rounding, and Diving \cite{berthold2014primal}. We turn on the aggressive mode in SCIP to focus on improving the primal bound instead of proving optimality, as the goal of this work is to develop efficient primal heuristics. Note that both Undercover and RENS could be automatically triggered as a subsidiary method during BnB \cite{BestuzhevaEtal2021OO}. We provide two SCIP baselines: 
    \begin{itemize}
        \item \textbf{SCIP-LP-form}. SCIP linearizes products of binary variables by default \cite{BestuzhevaEtal2021OO}, allowing the solver to employ widely studied MILP primal heuristics in solving MBQPs, as discussed in the background section.
        \item \textbf{SCIP-NLP-form}. In addition, we also run SCIP using the original MBQP formulation. 
    \end{itemize}
    
    \item \textbf{RENS}. We also compare with RENS as a standalone heuristic, as RENS could be applied along with other heuristics in the SCIP baselines or not triggered by SCIP. We use the official implementation of RENS in SCIP that contains additional advanced algorithmic components. To run standalone RENS in SCIP, we apply RENS at the root node and disable all other primal heuristics. Additionally, we remove restrictions for RENS, such as the minimum ratio of integer variables to be fixed, as SCIP uses such criteria to disable RENS when the subproblem is too large, when RENS is treated as a subsidiary method inside BnB \cite{berthold2014rens}. 
    
    Moreover, the pure RENS heuristic terminates after solving the relaxation and the sub-problem, which could happen before the time limit in our experiments. To give the RENS baselines more advantage, we allow SCIP to improve the solution returned by RENS through BnB until reaching the time limit (i.e., use RENS solution as an incumbent and continue branching after the root node).  
        \begin{itemize}
            \item \textbf{RENS\textsubscript{LP}+}. By default, SCIP in RENS uses the LP relaxation solution for fixing values. 
            \item \textbf{RENS\textsubscript{NLP}+}. We also add another RENS baseline that uses the NLP relaxation solution for fixing values.
        \end{itemize}
    We use the ``+" signs to highlight that we allow SCIP to further improve the solution returned by RENS through branching, which is a stronger baseline than pure RENS.

    \item \textbf{Undercover}. Similarly, we compare with Undercover as a standalone heuristic. We use the official full-fledged Undercover implemented in SCIP that contains additional features such as fix-and-propagate, backtracking, and postprocessing \cite{berthold2014undercover}. 
    \begin{itemize}            
            \item \textbf{Undercov\textsubscript{LP}+}. By default, Undercover in SCIP prioritizes using the LP relaxation as the fixing values for the cover.
            \item \textbf{Undercov\textsubscript{NLP}+} . We also run Undercover with NLP relaxation as the prioritized fixing values.
        \end{itemize}
        Again, we use ``+" to highlight that we allow SCIP to improve over the solution returned by Undercover after Undercover terminates. 
\end{itemize}






\subsubsection{Proposed Methods}
We experiment with both LP and NLP relaxation as the guidance for rounding in our proposed methods, resulting in the following variants:
\begin{itemize}
            \item \textbf{Relax-Search\textsubscript{LP}} and \textbf{Cover-Relax-Search\textsubscript{LP}}. We use SCIP to obtain the LP reformulation of the MBQPs and the LP relaxation solution. The time cutoff for solving the relaxation is $T_r=20s$. 
         
            \item \textbf{Relax-Search\textsubscript{NLP} and \textbf{Cover-Relax-Search\textsubscript{NLP}}}. We set a  limit of $T_r=20s$ for the NLP relaxation. However, because the trace of feasible solutions for the NLP relaxation is sparse compared to LP relaxation, 
            we stop solving the relaxation after the first relaxation solution is found.
        \end{itemize}
For the Cover-Relax-Search variants, we give a time limit of $T_c=1$s for solving the cover problem in (\ref{eq:1}). In all variants of our proposed methods, we use the rounded LP/NLP relaxation as a warm start solution when solving the sub-MBQPs. In creating the sub-MBQPs, we started with $p=0.7$ on the synthetic benchmarks and performed a sensitivity analysis with $p \in \{0.5, 0.6, 0.8, 0.9, 1\}$. The results of the sensitivity analysis showed that while the best value of $p$ is different for each benchmark, our conclusion is valid for all $p \in \{0.5, 0.6, 0.7, 0.8, 0.9, 1\}$. Therefore, we used $p=0.7$ in all synthetic benchmarks and applied the same value of $p$ to the real-world wind farm layout optimization problem instances in the next section.

\subsubsection{Computational Setup}
For all methods, we set the time limit to 60s. We conduct our experiments on 2.5 GHz AMD EPYC 7502 CPUs with 256 GB RAM. The proposed methods are implemented in Python. We use SCIP (v8.0.1) \cite{BestuzhevaEtal2021OO} to solve the cover problem, the LP and NLP relaxation problems, and the sub-MBQPs in our proposed methods.


\begin{table*}[h]
\centering

\scalebox{0.77}{
\begin{tabular}{clrrrrrr}
\Xhline{3\arrayrulewidth}
\\[-0.8em]
 &  &    & \textbf{Primal gap}  &    &   & \textbf{Primal integral}  & \\
\cmidrule(r){3-5}                                     \cmidrule(r){6-8} 
\\[-0.85em]
 & \textbf{Method}    &  CBQP  & CQKP  & QMKP   &  CBQP  & CQKP  & QMKP \\[-0.85em] \\ \Xhline{2\arrayrulewidth} \\[-0.85em]
\multirow{2}{*}{\begin{tabular}[c]{@{}c@{}}SCIP\\ Baselines\end{tabular}}       
 & SCIP-LP-form   & 0.94 (0\%)     & 0.61 (0\%)      & 0.41 (0\%)      & 57.82 (0\%)     & 41.64 (0\%)      & 29.15 (0\%)    \\[-0.85em] \\ \cline{2-8} \\[-0.85em]
 & SCIP-NLP-form       & 0.84 (10.64\%) & 0.77 (-26.23\%) & 0.77 (-87.8\%)  & 52.75 (8.77\%)  & 48.52 (-16.52\%) & 49.61 (-70.19\%) \\[-0.85em] \\ \Xhline{2\arrayrulewidth} \\[-0.85em]

 \multirow{2}{*}{\begin{tabular}[c]{@{}c@{}}RENS\\ Baselines\end{tabular}}       
  & RENS\textsubscript{LP}+   & 1 (-6.38\%)    & 0.99 (-62.3\%)  & 0.75 (-82.93\%) & 60 (-3.77\%)    & 59.85 (-43.73\%) & 53.74 (-84.36\%) \\[-0.85em] \\ \cline{2-8} \\[-0.85em]
 & RENS\textsubscript{NLP}+    & 1 (-6.38\%)    & 0.99 (-62.3\%)  & 0.75 (-82.93\%) & 60 (-3.77\%)    & 59.86 (-43.76\%) & 53.75 (-84.39\%) \\[-0.85em] \\ \Xhline{2\arrayrulewidth} \\[-0.85em]

  \multirow{2}{*}{\begin{tabular}[c]{@{}c@{}}\textbf{Proposed}\\ \textbf{Relax-Search}\end{tabular}}       
 & Relax-Search\textsubscript{LP} & 0.84 (10.64\%) & 0.14 (77.05\%)  & 0.24 (41.46\%)  & 54.3 (6.09\%)   & 21.9 (47.41\%)   & 21.34 (26.79\%) \\[-0.85em] \\ \cline{2-8} \\[-0.85em]
 & Relax-Search\textsubscript{NLP}    & \textbf{0 (100\%)}      &\textbf{ 0.07 (88.52\%) } & \textbf{0.01 (97.56\%)}  & \textbf{2.59 (95.52\%)}  & \textbf{10.73 (74.23\%) } & \textbf{9.45 (67.58\%)}  \\[-0.85em] \\ \Xhline{2\arrayrulewidth} \\[-0.85em]
 
\multirow{2}{*}{\begin{tabular}[c]{@{}c@{}}Undercov\\ Baselines\end{tabular}}       
  & Undercov\textsubscript{LP}+    & 0.94 (0\%)     & 0.61 (0\%)      & 0.44 (-7.32\%)  & 57.83 (-0.02\%) & 42.82 (-2.83\%)  & 30.52 (-4.7\%) \\[-0.85em] \\ \cline{2-8} \\[-0.85em]
 & Undercov\textsubscript{NLP}+  & 0.94 (0\%)     & 0.86 (-40.98\%) & 0.64 (-56.1\%)  & 58.12 (-0.52\%) & 53.68 (-28.91\%) & 41.19 (-41.3\%)  \\[-0.85em] \\ \Xhline{2\arrayrulewidth} \\[-0.85em]

 \multirow{2}{*}{\begin{tabular}[c]{@{}c@{}}\textbf{Proposed}\\ \textbf{Cover-Relax-Search}\end{tabular}}       
  & Cover-Relax-Search\textsubscript{LP}    & 0.83 (11.7\%)  & 0.16 (73.77\%)  & 0.23 (43.9\%)   & 54.79 (5.24\%)  & 23.1 (44.52\%)   & 23.12 (20.69\%)\\[-0.85em] \\ \cline{2-8} \\[-0.85em]
 & Cover-Relax-Search\textsubscript{NLP}  & 0.01 (98.94\%) & 0.07 (88.52\%)  & 0.02 (95.12\%)  & 3.96 (93.15\%)  & 11.74 (71.81\%)  & 12.09 (58.52\%) \\[-0.85em] \\ \Xhline{3\arrayrulewidth}
\end{tabular}
}
\caption{
\textbf{Results on standard MBQP benchmarks-small.} Primal gap (lower the better), Primal integral (lower the better), and the \% improvement compared to SCIP-LP-form (larger improvement is better) at 60s time cutoff, averaged over 100 test instances for each benchmark. The SCIP baselines can employ multiple primal heuristics, including RENS and Undercover.
}\label{tab:res_standard_small}
\end{table*}

\begin{table*}[h]
\centering

\scalebox{0.77}{
\begin{tabular}{clrrrrrr}
\Xhline{3\arrayrulewidth}
\\[-0.8em]
 &    &   & \textbf{Primal gap}  &    &   & \textbf{Primal integral}  & \\
\cmidrule(r){3-5}                                     \cmidrule(r){6-8}   
\\[-0.85em]
 & \textbf{Method}    &  CBQP  & CQKP  & QMKP   &  CBQP  & CQKP  & QMKP \\[-0.85em] \\ \Xhline{2\arrayrulewidth} \\[-0.85em]
 
\multirow{2}{*}{\begin{tabular}[c]{@{}c@{}}SCIP\\ Baselines\end{tabular}}       
 & SCIP-LP-form & 1 (0\%)     & 0.97 (0\%)     & 0.54 (0\%)      & 60 (0\%)        & 58.15 (0\%)     & 37.33 (0\%)  \\[-0.85em] \\ \cline{2-8} \\[-0.85em]
 & SCIP-NLP-form     & 1 (0\%)     & 0.91 (6.19\%)  & 1 (-85.19\%)    & 60 (0\%)        & 58.68 (-0.91\%) & 60 (-60.73\%) \\[-0.85em] \\ \hline \\[-0.85em]

 \multirow{2}{*}{\begin{tabular}[c]{@{}c@{}}RENS\\ Baselines\end{tabular}}       
  & RENS\textsubscript{LP}+     & 1 (0\%)     & 1 (-3.09\%)    & 0.95 (-75.93\%) & 60 (0\%)        & 60 (-3.18\%)    & 59.52 (-59.44\%) \\[-0.85em] \\ \cline{2-8} \\[-0.85em]
 & RENS\textsubscript{NLP}+   & 1 (0\%)     & 1 (-3.09\%)    & 0.95 (-75.93\%) & 60 (0\%)        & 60 (-3.18\%)    & 59.52 (-59.44\%)  \\[-0.85em] \\ \hline \\[-0.85em]

  \multirow{2}{*}{\begin{tabular}[c]{@{}c@{}}\textbf{Proposed}\\ \textbf{Relax-Search}\end{tabular}}       
 & Relax-Search\textsubscript{LP}  & 0.93 (7\%)  & 0.5 (48.45\%)  & 0.26 (51.85\%)  & 58.13 (3.12\%)  & 39.76 (31.63\%) & 37.05 (0.75\%)  \\[-0.85em] \\ \cline{2-8} \\[-0.85em]
 & Relax-Search\textsubscript{NLP}     & \textbf{0.01 (99\%)} & \textbf{0.02 (97.94\%)} & \textbf{0.05 (90.74\%)}  & \textbf{15.2 (74.67\%)}  & \textbf{22.58 (61.17\%)} & \textbf{26.7 (28.48\%)} \\[-0.85em] \\ \hline \\[-0.85em]
 
\multirow{2}{*}{\begin{tabular}[c]{@{}c@{}}Undercover\\ Baselines\end{tabular}}       
  & Undercov\textsubscript{LP}+     & 1 (0\%)     & 1 (-3.09\%)    & 1 (-85.19\%)    & 60 (0\%)        & 60 (-3.18\%)    & 60 (-60.73\%)  \\[-0.85em] \\ \cline{2-8} \\[-0.85em]
 & Undercov\textsubscript{NLP}+  & 1 (0\%)     & 1 (-3.09\%)    & 1 (-85.19\%)    & 60 (0\%)        & 60 (-3.18\%)    & 60 (-60.73\%)   \\[-0.85em] \\ \hline \\[-0.85em]

 \multirow{2}{*}{\begin{tabular}[c]{@{}c@{}}\textbf{Proposed}\\ \textbf{Cover-Relax-Search}\end{tabular}}       

  & Cover-Relax-Search\textsubscript{LP}+  & 0.94 (6\%)  & 0.62 (36.08\%) & 0.49 (9.26\%)   & 58.82 (1.97\%)  & 43.37 (25.42\%) & 40.74 (-9.13\%) \\[-0.85em] \\ \cline{2-8} \\[-0.85em]
 & Cover-Relax-Search\textsubscript{NLP}+  & \textbf{0.01 (99\%)} & 0.07 (92.78\%) & 0.13 (75.93\%)  & 17.79 (70.35\%) & 24.98 (57.04\%) & 29.1 (22.05\%) \\[-0.85em] \\ \Xhline{3\arrayrulewidth}
\end{tabular}
}
\caption{\textbf{Results on standard MBQP benchmarks-large.} Primal gap (lower the better), Primal integral (lower the better), and the \% improvement compared to SCIP-LP-form (larger improvement is better) at 60s time cutoff, averaged over 100 test instances for each benchmark. The SCIP baselines can employ multiple primal heuristics, including RENS and Undercover.}\label{tab:res_standard-large}
\end{table*}

\subsection{Results and Discussion}

The results for the small and large generated MBQP benchmarks are shown in Table \ref{tab:res_standard_small} and Table \ref{tab:res_standard-large}, respectively. Relax-Search\textsubscript{NLP} performs the best in all small and large benchmarks, with Cover-Relax-Search\textsubscript{NLP} being the second-best method. We report the improvement in primal gap and primal integral compared to the SCIP baseline with LP formulation (default setting) for all the methods. 

The standalone RENS and Undercover baselines (even when we allow SCIP to improve the resulting solution through BnB) fail to outperform SCIP-LP-form in all the tested benchmarks. This is expected, as SCIP employs a rich set of primal heuristics in addition to RENS and Undercover. Relax-Search\textsubscript{NLP} reduces the primal gap by $88.52-100\%$ and reduces the primal integral by $28.48-95.52\%$ compared to SCIP-LP-form in different benchmarks. Cover-Relax-Search\textsubscript{NLP} also significantly reduces the primal gap and primal integral compared to SCIP-LP-form. This means that our proposed methods find better solutions at the time cutoff and produce high-quality solutions at a faster rate.
\begin{figure*}[h]
\centering
\begin{subfigure}[b]
{0.49\textwidth}
\centering
   \includegraphics[width=1\linewidth]{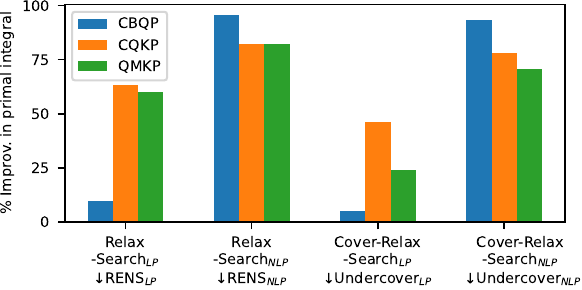}
   \caption{small}
   \label{fig:Ng1} 
\end{subfigure}
\begin{subfigure}[b]{0.49\textwidth}
\centering
   \includegraphics[width=1\linewidth]{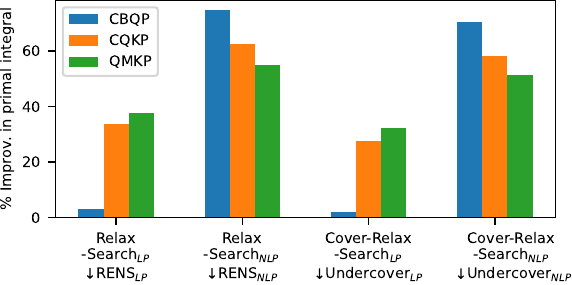}
   \caption{large}
   \label{fig:Ng2}
\end{subfigure}
\caption{\textbf{Benefits of using suboptimal relaxation as guidance for rounding and variable fixing on standard MBQP benchmarks.} \% improvement in primal integral at 60s using the proposed methods when compared to their most similar baseline methods.}
\end{figure*}
\paragraph{Benefits of Suboptimal Relaxation.}
The weak performance of RENS and Undercover on the tested MBQPs can likely be attributed to the high computational cost of solving the relaxation. In several cases, the primal gap at the time cutoff is 1 for the RENS and Undercover baselines, both with LP and NLP relaxation, as the basis for fixing values. A primal gap of 1 indicates that the method fails to find an improved solution beyond the initial trivial solution, given that SCIP uses 0 as the trivial starting solution in all our tested benchmarks. In these cases, RENS and Undercover spend a large portion of their time solving or improving the LP/NLP relaxation, leaving little or no time for solving the restricted subproblem. Fig. \ref{fig:Ng1} and \ref{fig:Ng2} further illustrate the benefits of using suboptimal relaxation as guidance for rounding and variable fixing. In these plots, we show the percentage improvement in primal integral using the proposed methods when compared to their most similar baseline methods (as noted in the x-axis labels). These results show that 
setting a time cutoff for solving the relaxation and using the suboptimal relaxation as a basis for creating the sub-MBQPs (i.e., the proposed methods) leads to better performance across all benchmarks, both when using LP and NLP relaxation. 



\paragraph{Benefits of NLP Relaxation.}
The NLP relaxation is less commonly used in rounding-based heuristics, as it is generally considered more computationally expensive than the LP relaxation. However, our results indicate that even a suboptimal NLP relaxation solution can yield better results than an LP relaxation when it is used as a guide for variable fixing. A possible explanation is that the LP reformulation relaxation for MBQPs is weak.


\section{Evaluation on Real-World Wind Farm Layout Optimization Problem}
We examine the efficacy of our proposed methods by applying them to a large-scale, real-world Wind Farm Layout Optimization Problem (WFLOP). Wind energy is harnessed in wind farms using wind turbines, which transform the kinetic energy of transient winds into electrical energy. 
Intuition suggests that to maximize the power produced by a wind farm, one should maximize the number of installed wind turbines. 
However, this is not the case due to a phenomenon called \emph{wake effects}.
Once installed, the layout of the turbines will have a significant impact on the wake effects and, therefore, the power production under uncertain wind conditions.
\cite{barthelmie2009modelling} estimate that in the case of large offshore wind farms, the average power lost due to turbine wake effects is between 10-20\% the annual power production.
Therefore, it is critically important to identify optimal wind turbine layouts that maximize power production by minimizing wind speed losses caused by wake effects. 
\subsection{Problem Description}

WFLOP is a class of design optimization problems concerned with the placement of $K$ wind turbines within a specified area, or in a finite set of discrete locations, to maximize total wind farm power production. 
This problem has been studied extensively using various formulations, including continuous nonlinear programming \cite{perez2013offshore} and mixed-integer models \cite{turner2014new}. In the following section, we extend the original deterministic model in \cite{turner2014new} to include an expectation objective over uncertain wind conditions.

\paragraph{MBQP Formulation}

Let $J$ be the set of candidate locations for turbine placement and $K$ be the number of turbines to be installed. The binary decision variable $y_j$ takes the value 1 if a wind turbine is installed at location $j$, and 0 otherwise. 
The set $J$ is a two-dimensional grid with a given resolution that represents the design area of the wind farm. 

The free stream wind speed is represented by parameter $U$ and the wind direction is $\theta$. Assume that we have a set of $M$ wind scenarios $\mathcal{M} = \{1,2,\ldots,M\}$ drawn from a joint probability distribution $p(U,\theta)$. Each scenario $m$ consists of a wind speed $U^{(m)}$ and wind direction $\theta^{(m)}$ with  probability $p^{(m)}$ such that $\sum_{m \in \mathcal{M}} p^{(m)}  = 1$. 
The pairwise wind speed deficit interactions $d_{ij}^{(m)}$ are a function of additional parameters, including wind direction $\theta^{(m)}$. 
These are precomputed for a given grid resolution and wind scenario set $\mathcal{M}$.

Given these modeling assumptions, the WFLOP MBQP model to minimize expected wind speed losses due to wake interactions across all turbines in the wind farm is shown in Equations~(\ref{eqn:bqp-model-obj})--(\ref{eqn:bqp-model-constr}). 

\begin{alignat}{2}
\text{min} \quad & \sum_{m \in \mathcal{M}} p^{(m)}\, U^{(m)} \sum_{j \in J}\sum_{i \in J} \left(d_{ij}^{(m)}\right)^2\, y_i  y_j \label{eqn:bqp-model-obj}\\
\text{s.t.:} \quad & \sum_{j \in J} y_j = K \\
& y_j \in \{0,1\} \quad \forall j \in J\label{eqn:bqp-model-constr}
\end{alignat}

\begin{figure}[h]
\centering
\begin{subfigure}[b]{.4\linewidth}
\includegraphics[width={1\textwidth}]{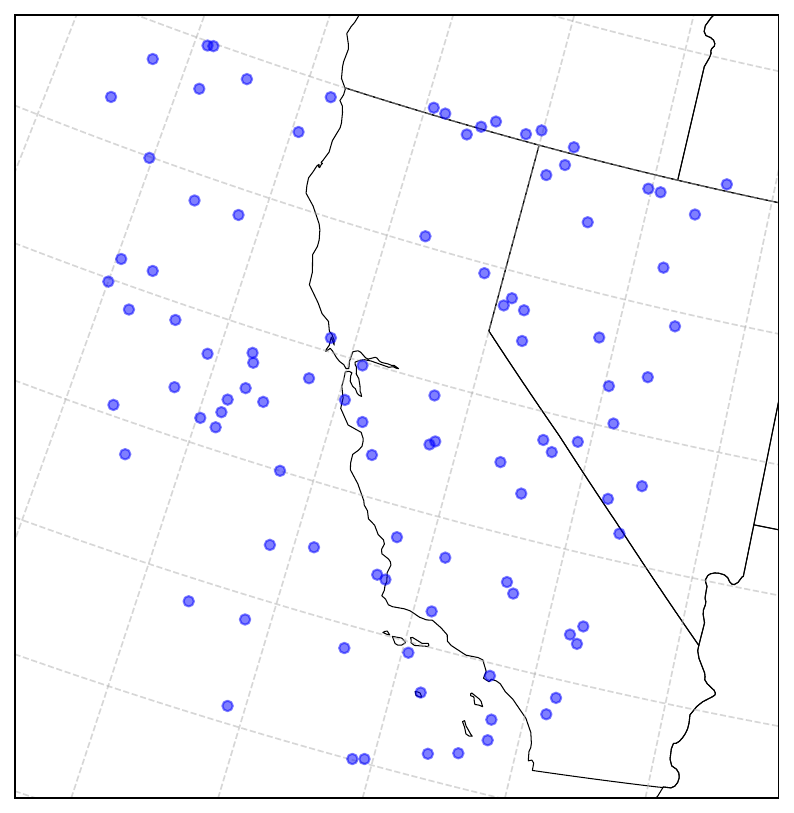}
\caption{}\label{fig:sites}

\end{subfigure}
\begin{subfigure}[b]{.56\linewidth}

\includegraphics[width={1\textwidth}]{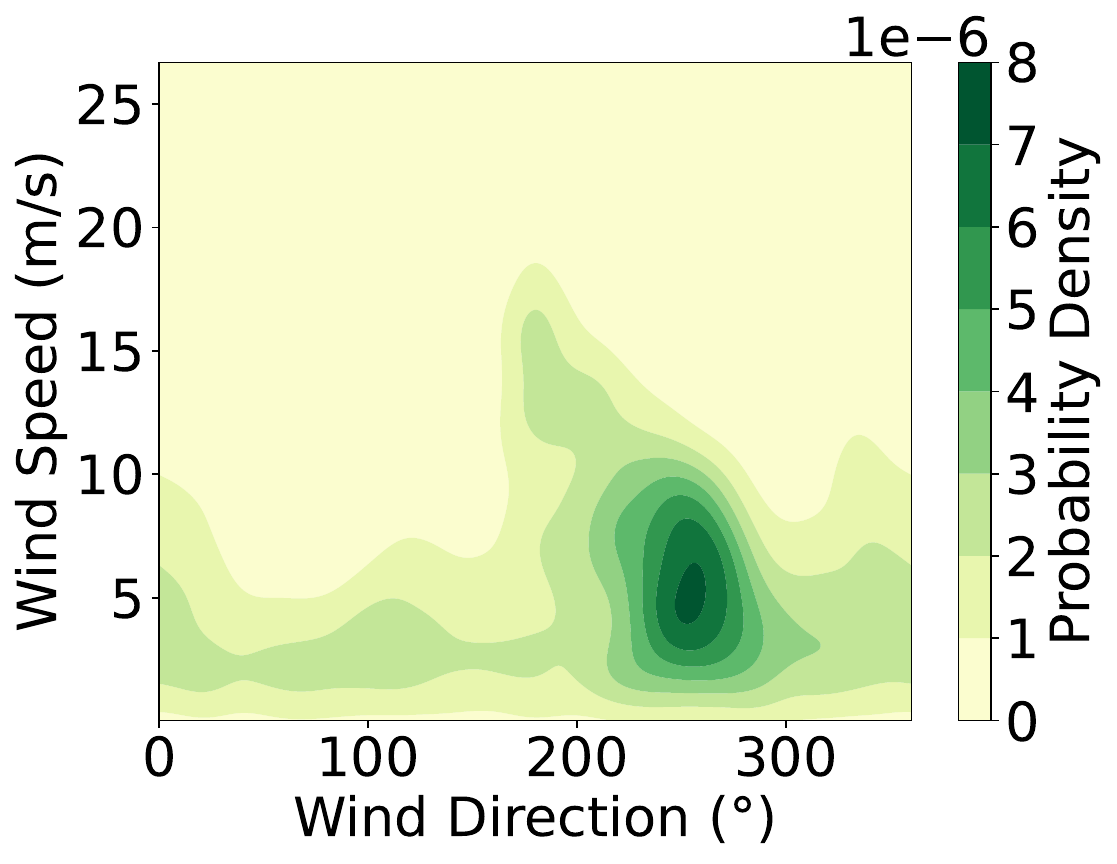}
\caption{}\label{fig:pdf}
\end{subfigure}
\caption{(a) 100 geographic sites selected from the NOW-23 California offshore wind data to instantiate the WFLOP MBQP instances (1 small and 1 large WFLOP instance at each dot). (b) Representative wind joint probability distribution kernel density estimate for $p(U,\theta)$ for a given site.}
\end{figure}

\paragraph{Data} 
%
%
We use wind data from the National Energy Technology Laboratory (NREL) 2023 National Offshore Wind data set (NOW-23) \cite{bodini20232023}, which contains annual temporal wind data simulated using the Weather Research and Forecasting model \cite{skamarock2019description} over several different offshore regions in the United States. For this study, we use hourly resolution data for year 2019 within the California geographic region. 100 locations in California (defined by longitude and latitude) were randomly selected from the simulated wind data. 
At each location, the wind speed and direction data at 100 meter height (e.g., approximate wind turbine height) were used to fit a joint probability distribution $p(U,\theta)$ using kernel density estimation.  
The sites used in our study are shown in Fig.~\ref{fig:sites}.
An example joint probability distribution of wind speed $U$ and direction $\theta$ for a given location in Fig.~\ref{fig:sites} is shown in Fig.~\ref{fig:pdf}.

\subsubsection{Setup}
We create two sets of WFLOP benchmarks, referred to as \emph{small} and \emph{large}, each containing 100 instances that correspond to the 100 locations selected in California (Fig.~\ref{fig:sites}). 
Each of these locations experiences different wind phenomena, leading to 100 unique benchmarks based on the underlying joint distribution of wind speed and direction, i.e., $p(U,\theta)$.
The instances in the small set are defined for a single wind scenario drawn randomly from the site-specific joint distribution, making the set $\mathcal{M}$ a singleton. 
The large instances are defined for 10 wind scenarios randomly drawn from the site-specific joint distribution, i.e., $\mathcal{M}=\{1,\ldots,M\}$ with $M=10$.
The average quadratic term density in Equation~(\ref{eqn:bqp-model-obj}) for small and large instances is 0.075 and 0.36, respectively. The same computational settings applied to the synthetic benchmarks were also used for the WFLOP instances.

\subsubsection{Results}
Relax-Search\textsubscript{NLP} results in the best primal gap and primal integral in both the small and large WFLOP benchmarks. Additionally, the proposed Cover-Relax-Search\textsubscript{NLP} also leads to the best primal gap in the small instances set.

\begin{table}[!tb]
\centering

\scalebox{0.78}{
\begin{tabular}{clrrrr}
\Xhline{3\arrayrulewidth}
\\[-0.8em]
 &  &  \multicolumn{2}{c}{\textbf{Small}}  &   \multicolumn{2}{c}{\textbf{Large}}   \\
\cmidrule(r){3-4}                                     \cmidrule(r){5-6} 
\\[-0.85em]
 & \textbf{Method}   & PG  & PI     & PG  & PI     \\[-0.85em] \\ \hline \\[-0.85em]
\multirow{2}{*}{\begin{tabular}[c]{@{}c@{}}SCIP\\ Baselines\end{tabular}}       
 & SCIP-LP-form    & 0.27                   & 24.1                   & 0.19 & 18.48  \\[-0.85em] \\ \cline{2-6} \\[-0.85em]
 & SCIP-NLP-form   & 0.44                    & 29.84                  &0.49 & 32.7   \\[-0.85em] \\ \hline \\[-0.85em]
 \multirow{2}{*}{\begin{tabular}[c]{@{}c@{}}RENS\\ Baselines\end{tabular}}       
  & RENS\textsubscript{LP}+     & 0.34                   & 30.06                  & 0.3  & 26.24 \\[-0.85em] \\ \cline{2-6} \\[-0.85em]
 & RENS\textsubscript{NLP}+   & 0.34                   & 30.07                  & 0.3  & 26.23 \\[-0.85em] \\ \hline \\[-0.85em]

 \multirow{2}{*}{\begin{tabular}[c]{@{}c@{}}\textbf{Proposed}\\ \textbf{Methods}\end{tabular}}       
 & Relax-Search\textsubscript{LP}  & 0.5                   & 31.16                  & 0.48 & 30.64  \\[-0.85em] \\ \cline{2-6} \\[-0.85em]
 & Relax-Search\textsubscript{NLP}    & \textbf{0.01}                      & \textbf{3.34}                   & \textbf{0.03} & \textbf{14.53} \\[-0.85em] \\ \hline \\[-0.85em]

\multirow{2}{*}{\begin{tabular}[c]{@{}c@{}}Undercover\\ Baselines\end{tabular}}       
  & Undercov\textsubscript{LP}+   & 0.59                   & 36.49                 & 0.55 & 34.08  \\[-0.85em] \\ \cline{2-6} \\[-0.85em]
 & Undercov\textsubscript{NLP}+   & 0.6                   & 60                  & 0.55 & 59.99\\[-0.85em] \\ \hline \\[-0.85em]

 \multirow{2}{*}{\begin{tabular}[c]{@{}c@{}}\textbf{Proposed}\\ \textbf{Methods}\end{tabular}}       
  & Cover-Relax-Search\textsubscript{LP}  & 0.37                   & 24.39                  & 0.48 & 30.77   \\[-0.85em] \\ \cline{2-6} \\[-0.85em]
 & Cover-Relax-Search\textsubscript{NLP}  & \textbf{0.01}                      & 4.19                      & 0.04 & 15.96 \\[-0.85em] \\ \Xhline{3\arrayrulewidth}
\end{tabular}
}
\caption{\textbf{Results on real-world wind farm layout problem.} Primal Gap (PG) (lower the better) and Primal Integral (PI) (lower the better) averaged over
100 test instances.}\label{tab:res_wflop}
\end{table}
\begin{figure*}[h]
\centering
\begin{subfigure}{.33\textwidth}
\centering
\includegraphics[width={0.95\textwidth}]{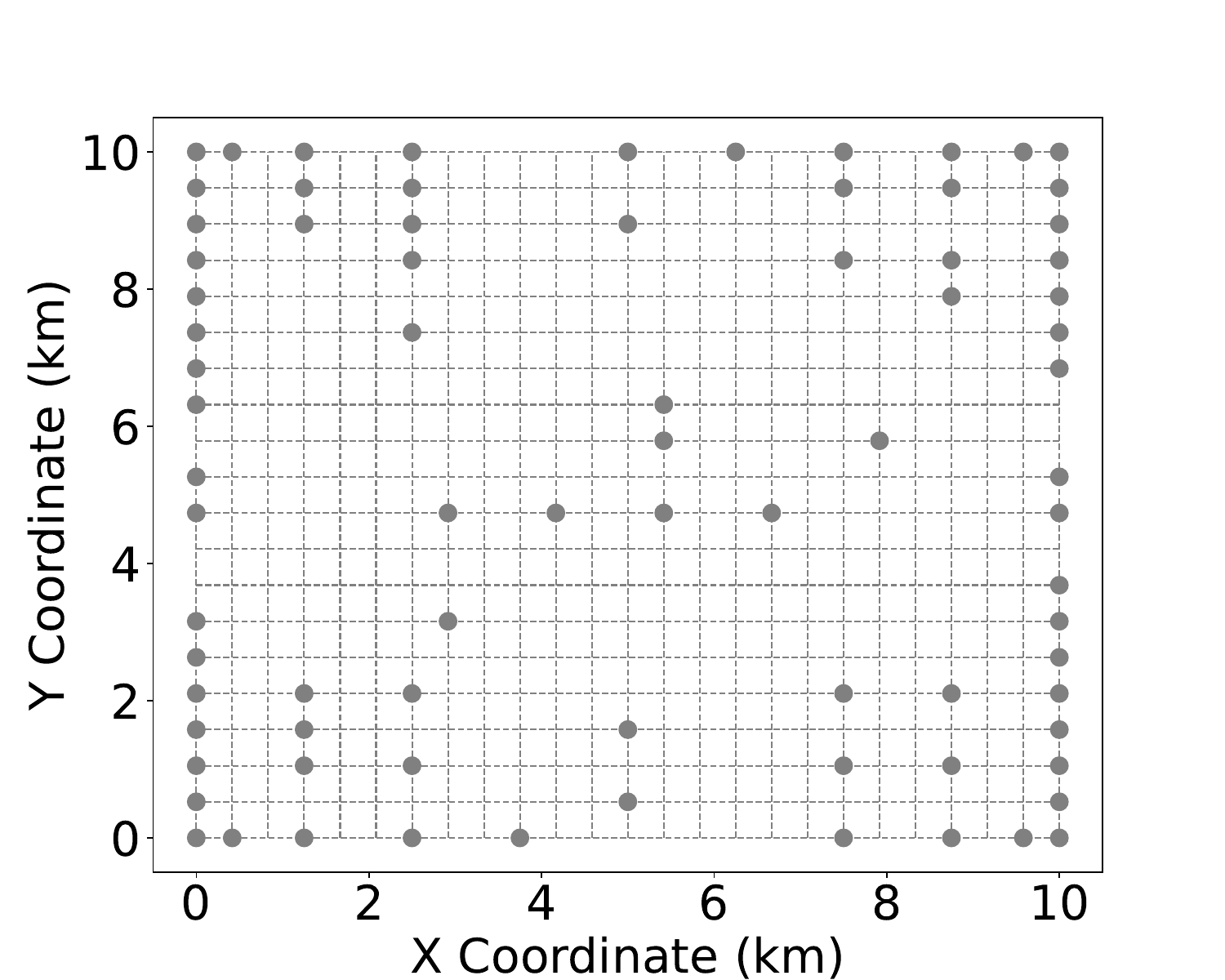}
\caption{Relax-Search\textsubscript{NLP} (ours)}\label{fig:layout-relax}
\end{subfigure}
\begin{subfigure}{.33\textwidth}
\centering
\includegraphics[width={0.95\textwidth}]{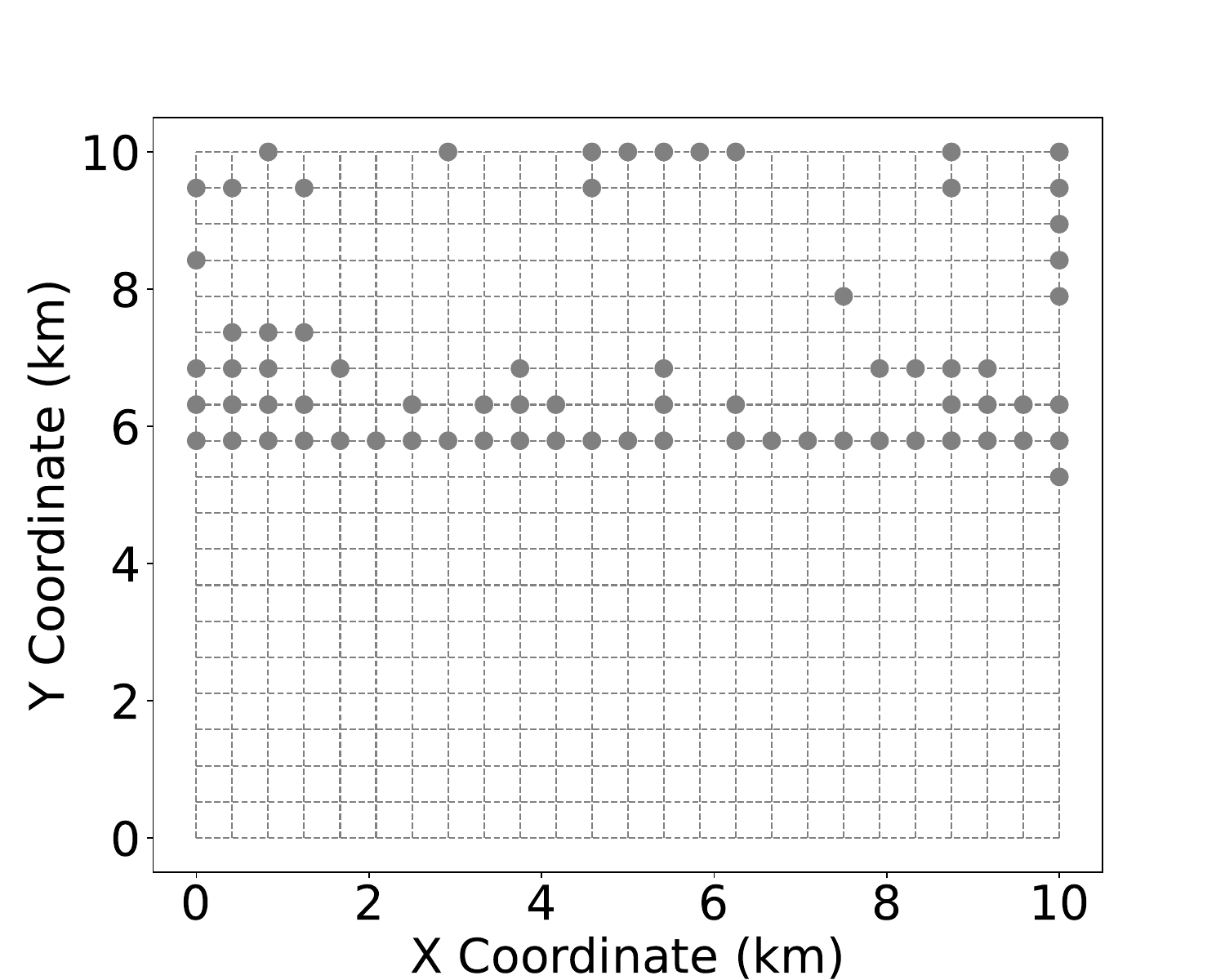}
\caption{RENS\textsubscript{NLP}+}\label{fig:layout-rens}
\end{subfigure}
\begin{subfigure}{.33\textwidth}
\centering
\includegraphics[width={0.95\textwidth}]{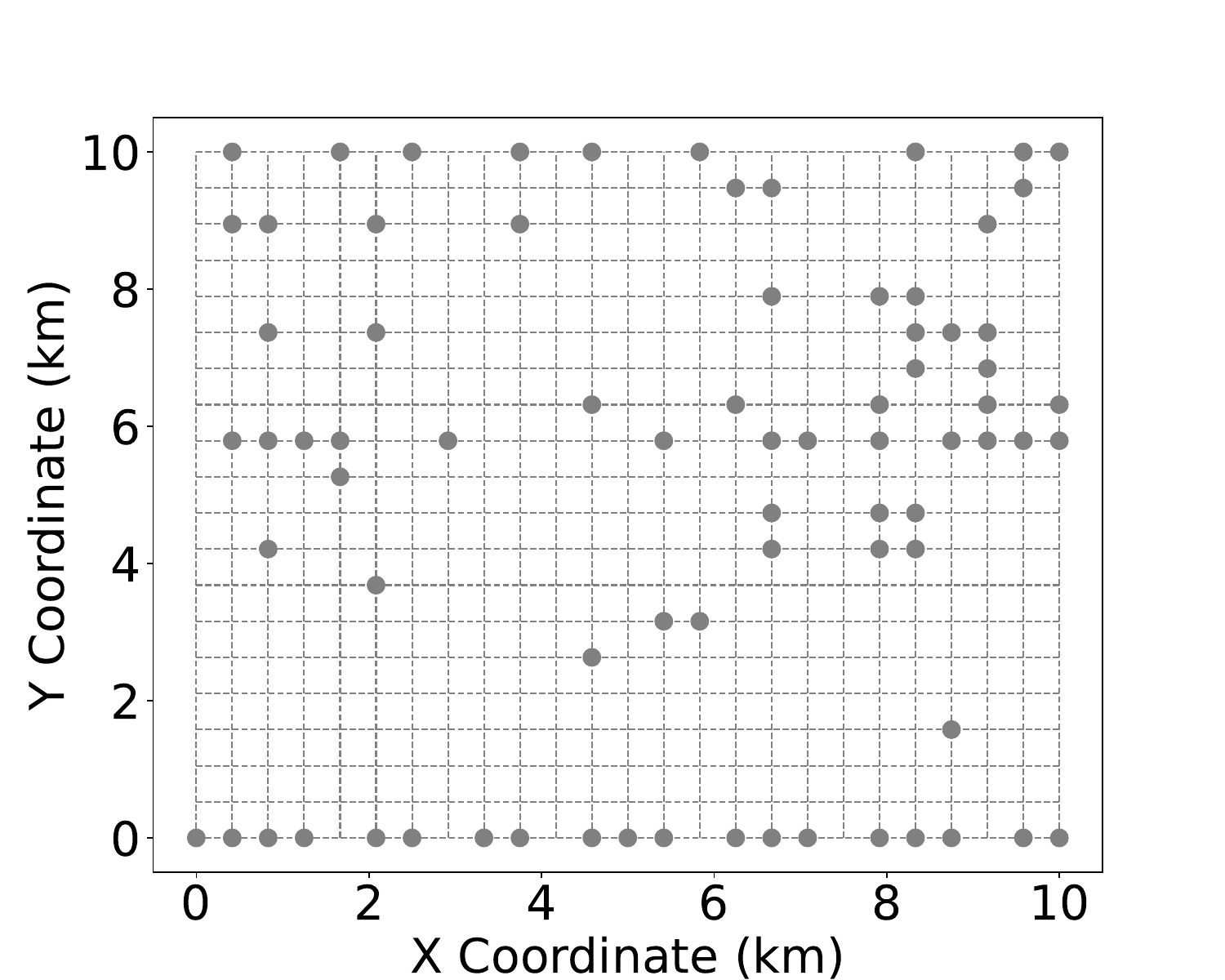}
\caption{SCIP-LP-form}\label{fig:layout-bnb}
\end{subfigure}
\caption{Final turbine layout designs returned by each respective algorithm for a given large WFLOP instance. Grids represent candidate locations while scatter points represent allocated turbine locations.}\label{fig:layouts}
\end{figure*}
In Fig.~\ref{fig:layouts}, we plot an example of turbine layout solutions identified by Relax-Search\textsubscript{NLP}, RENS\textsubscript{NLP}+, and SCIP-LP-form on a large WFLOP instance.
These solutions represent the best overall objective value and primal integral, the best-performing RENS baseline, and the best-performing SCIP baseline from the results in Table~\ref{tab:res_wflop}. One notable difference between these turbine layout solutions is the relative sparsity of turbine placement. 
In Fig.~\ref{fig:layout-relax}, the turbines are placed in a more diffuse pattern than in \ref{fig:layout-rens} and \ref{fig:layout-bnb}. 
The latter two designs contain more closely clustered regions of directly adjacent wind turbines. 
Additionally, the layout solution in Fig.~\ref{fig:layout-relax} appears to prioritize placing turbines along the perimeter, maximizing inter-turbine spacing.
The relative improvement of the Relax-Search\textsubscript{NLP} design over the SCIP-LP-form design results in 8\% reduction in wind speed losses, corresponding to a 32\% increase in expected power production.

This dramatic improvement can be explained through the underlying wake interaction model, wherein wake effects dissipate as a function of inter-turbine distance \cite{turner2014new}.
Thus, it is generally advantageous to maximize spacing between turbines to mitigate wake effects. 
This explains why the sparser design identified by Relax-Search\textsubscript{NLP} significantly outperforms the next best baseline solution. 
Furthermore, our results show that identifying high-quality turbine layout designs is not straightforward in the case of uncertain wind conditions.
This is because inter-turbine distances are computed in the direction of the prevailing winds and must be averaged over multiple scenarios. 

\section{Related Work}\label{sec:related}

To our knowledge, there has been limited work in the existing literature on primal heuristics specifically for solving MBQPs; many are developed for general MINLPs. Primal heuristics for MINLPs are often adaptations of MILP heuristics, including the Relaxation Induced Neighborhood Search \cite{berthold2014primal}, Relaxation Enforced Neighborhood Search \cite{berthold2014rens}, Feasibility Pump \cite{bonami2009feasibility}, and Local Branching \cite{berthold2014primal}. Some studies have proposed algorithms for specific classes of MBQPs. \cite{gomez2024real} introduced a method for obtaining real-time solutions to MBQPs with banded matrices and indicator variables by using decision diagrams and solving a shortest-path problem. \cite{takapoui2020simple} developed a heuristic to identify approximate solutions to MBQPs with a convex objective function based on the alternating direction method of multipliers. 

\section{Conclusion and Discussion}\label{sec:conclusions}
In this work, we develop primal heuristics to efficiently identify high-quality solutions for large-scale MBQPs. We propose \emph{Relax-Search} and \emph{Cover-Relax-Search}, which extend the RENS and Undercover heuristics, respectively. Our methods use the suboptimal relaxation solutions to the original MBQP as guidance, selectively fix a subset of binary variables controlled by a ratio, and solve a restricted sub-MBQP to improve upon the suboptimal rounding. We evaluate our methods on both standard and real-world MBQP instances. To construct real-world MBQP benchmarks, we propose a formulation of the wind farm layout optimization problem that captures uncertain wind conditions and use wind distribution data in California to instantiate the models. Experimental results show that the proposed Relax-Search\textsubscript{NLP} method outperforms RENS, Undercover, and SCIP with additional primal heuristics integrated, achieving the lowest average primal gap and primal integral in all tested benchmarks. Furthermore, our analysis shows that while NLP relaxation is less commonly used in rounding-based primal heuristics due to the high 
computational cost, using suboptimal NLP relaxation can significantly enhance the performance. For future work, we plan to extend our methods to Mixed-Integer Quadratically Constrained Programs (MIQCPs), which is a broader class of problems that captures challenging optimization problems arising in many scientific domains.

\section{Acknowledgments}
This work was done during Weimin Huang’s internship at the Pacific Northwest National Laboratory (PNNL) and at the University of Southern California. PNNL is a multi-program national laboratory operated by Battelle Memorial Institute for the U.S. Department of Energy (DOE) under Contract No. DE-AC05-76RL0-1830. The research is partially supported by the National Science Foundation (NSF) grant 2112533: “NSF Artificial Intelligence (AI) Research Institute for Advances in Optimization (AI4OPT)”,  the U.S. Department of Energy, Office of Science Energy Earthshot Initiative, as part of the Addressing Challenges in Energy: Floating Wind in a Changing Climate Energy Earthshot Research Center at PNNL, and the Ralph S. O’Connor Sustainable Energy Institute (ROSEI) at Johns Hopkins University. 
\bibliography{aaai25}

\end{document}